\documentstyle[amssymb,amsfonts]{amsart}

\def\normo#1{\left\|#1\right\|}

\def\aabs#1{\left|#1\right|}
\def\brk#1{\left(#1\right)}

\def\norm#1{\|#1\|}
\def\jb#1{\langle#1\rangle}

\def\wh#1{\widehat{#1}}

\numberwithin{equation}{section}

\newcommand{\N}{{\mathbb N}}

\newcommand{\R}{{\mathbb R}}

\newcommand{\Z}{{\mathbb Z}}
\newcommand{\ft}{{\mathcal{F}}}

\newcommand{\les}{{\lesssim}}
\newcommand{\ges}{{\gtrsim}}

\newcommand{\Sch}{{\mathcal{S}}}

\theoremstyle{plain}
  \newtheorem{theorem}[subsection]{Theorem}
  \newtheorem{proposition}[subsection]{Proposition}
  \newtheorem{lemma}[subsection]{Lemma}
  \newtheorem{corollary}[subsection]{Corollary}

\theoremstyle{remark}
  \newtheorem{remark}[subsection]{Remark}

\theoremstyle{definition}

\newenvironment{proof}{\noindent {\bf Proof.} }{\endprf\par}
\def \endprf{\hfill  {\vrule height6pt width6pt depth0pt}\medskip}

\begin{document}

\title[Well-posedness for the Schr\"odinger-Korteweg-de Vries system]{On the  Well-posedness of the Schr\"odinger-Korteweg-de Vries system}
\author{Zihua Guo, Yuzhao Wang}
\address{LMAM, School of Mathematical Sciences, Peking University, Beijing
100871, China}

\email{zihuaguo, wangyuzhao@@math.pku.edu.cn}

\begin{abstract}
We prove that the Cauchy problem for the Schr\"odinger-Korteweg-de
Vries system is locally well-posed for the initial data belonging to
the Sovolev spaces $L^2(\R)\times H^{-{3/4}}(\R)$. The new
ingredient is that we use the $\bar{F}^s$ type space, introduced by
the first author in \cite{G}, to deal with the KdV part of the
system and the coupling terms. In order to overcome the difficulty
caused by the lack of scaling invariance, we prove uniform estimates
for the multiplier. This result improves the previous one by Corcho
and Linares \cite{cl07}.
\end{abstract}


\maketitle

\section{Introduction}\label{sect:intro_main}

In this paper, we consider the Cauchy problem for the
Schr\"odinger-Korteweg-de Vries (NLS-KdV) system
\begin{align}\label{S-KdV}
\begin{cases}
i\partial_tu +\partial^2_x u =\alpha uv+ \beta |u|^2u, \qquad t,x\in
\R,\\
\partial_t v +\partial_{x}^{3}v+\frac{1}{2}\partial_x(v^2)=\gamma\partial_x(|u|^2),\\
u(x,0)=u_{0}(x), v(x,0)=v_{0}(x),
\end{cases}
\end{align}
where $u(x,t)$ is a complex-valued function, $v(x,t)$ is a
real-valued function and $\alpha$, $\beta$, $\gamma$ are real
constants, $(u_0,v_0)$ are given initial data belonging to
$H^{s_1}\times H^{s_2}$. Our main motivation of this paper is
inspired by the work of Corcho and Linares \cite{cl07}, and the work
of the first named author \cite{G}.

The system \eqref{S-KdV} is an important model in fluid mechanics
and plasma physics that governs the interactions between short-wave
and long wave. The case $\beta=0$ appears in the study of resonant
interaction between short and long capillary-gravity waves on water
of uniform finite depth, in plasma physics and in a diatomic lattice
system (See \cite{cl07} and reference therein for more
introduction).

Before we state our main results, we recall first the early results
on this system. M. Tsutsumi \cite{t93} obtained global
well-posedness (GWP) for data $(u_0, v_0) \in H^{s+1/2}(\R)\times
H^s(\R)$ with $s \in Z_+$. In the resonant case $(\beta = 0)$ Guo
and Miao \cite{gm99} proved GWP in the natural energy space
$H^{s}(\R)\times H^s(\R)$ with $s \in Z_+$. Bekiranov, Ogawa and
Ponce \cite{bop97} proved local well-posedness (LWP) in
$H^{s}(\R)\times H^{s-1/2}(\R)$ for $s\geq 0$, and this result was
improved to $L^2\times H^{-3/4+}$ by Corcho and Linares \cite{cl07}
which seems to be sharp except $L^2\times H^{-3/4}$ in view of the
results for the two single equations in \eqref{S-KdV} (See
\cite{KPVDMJ01}, \cite{cct03}, \cite{KPVJAMS96}). Pecher \cite{p}
obtained GWP in $H^s\times H^s$ for $s>3/5$ ($\beta=0$) and $s>2/3$
($\beta\ne 0$) by using the ideas of I-method \cite{I-method}. Some
generalized interaction equations were considered in \cite{bop98}.
In this paper, we prove the following results:
\begin{theorem}\label{thmlwp} Assume $u_0\in L^2$, $v_0 \in H^{-3/4}$. Then

(a) Existence. There exist $T=T(\norm{u_0}_{L^2},
\norm{v_0}_{H^{-3/4}})>0$ and a solution $u$ to the Cauchy problem
\eqref{S-KdV} satisfying
\[u\in X^{0,1/2+}_{\tau=-\xi^2}(T)\subset C([-T,T]:L^2),\, v\in \bar{F}_{\tau=\xi^3}^{-3/4}(T) \subset C([-T,T]:H^{-3/4}).\]

(b) Uniqueness. The solution mapping $S_T:(u_0,v_0)\rightarrow
(u,v)$ is the unique extension of the classical solution
$(H^\infty,H^\infty)\rightarrow C([-T,T]:H^\infty\times H^\infty)$.

(c) Lipschitz continuity. For any $R>0$, the mapping
$(u_0,v_0)\rightarrow (u,v)$ is Lipschitz continuous from
$\{(u_0,v_0)\in L^2\times
H^{-3/4}:\norm{u_0}_{L^2}+\norm{u_0}_{H^{-3/4}}<R\}$ to
$C([-T,T]:L^2\times H^{-3/4})$.
\end{theorem}

We describe briefly our ideas in proving Theorem \ref{thmlwp}. We
also use the scheme as in \cite{cl07} which is the same spirit as
the one by Ginibre, Y. Tsutsumi and Velo \cite{gtv97} for the
Zakharov system. The basic idea is that for the second equation in
\eqref{S-KdV} we use the $\bar{F}^s$ space that was used by the
first named author \cite{G} for the KdV equation, but there is an
essential difficulty. For the KdV equation, one can assume the
initial data $v_0$ has a small norm by using the scaling transform.
However, for the NLS-KdV system we don't have such an invariant
scaling transform. In order to deal with large initial data, we
overcome this difficulty by the following way. We observe that the
single nonlinear Schr\"odinger equation with cubic term ($|u|^2u$)
\[i\partial_tu +\partial^2_x u =\beta |u|^2u\]
is $L^2$-subcritical. On the other-hand, we also see from
\cite{cl07} that to control the coupling term $uv$ one need less
regularity than $H^{-3/4}$ of $v$. Then we expect that the first
equation can be handled without scaling. Thus we scale the system
\eqref{S-KdV} according to the second equation as following: if
$(u,v)$ solve the system \eqref{S-KdV} with initial data
$(u_0,v_0)$, then we see
\begin{align}\label{scaling}
u_{\lambda}(t,x)&=\lambda^2 u(\lambda^3 t,\lambda x),\quad
\phi_1(x)=\lambda^2 u_0(\lambda x),\\
v_{\lambda}(t,x)&=\lambda^2 v(\lambda^3 t,\lambda x),\quad
\phi_2(x)=\lambda^2 v_0(\lambda x),
\end{align}
satisfy the following system
\begin{align}\label{S-KdV1}
\begin{cases}
i\partial_tu+\lambda\partial^2_x u =\lambda uv+ \lambda^{-1} |u|^2u,
\quad t,x\in \R,
\\
\partial_t v +\partial_{x}^{3}v+\frac{1}{2}\partial_x(v^2)=
\partial_x(|u|^2),
\\
u(0,x)=\phi_1(x), v(0,x)=\phi_2(x).
\end{cases}
\end{align}
It is easy to see that by taking $0<\lambda\ll 1$, we have
$\norm{\phi_1}_{L^2}=\norm{\lambda^2 u_0(\lambda
x)}_{L^2}=\lambda^{3/2}\norm{u_0}_2\leq R$ and
$\norm{\phi_2}_{H^{-3/4}}=\norm{\lambda^2 v_0(\lambda
x)}_{H^{-3/4}}\leq 2\lambda^{3/4}\norm{v_0}_{H^{-3/4}}=\epsilon_0\ll
1$. Therefore, it reduces to study the system \eqref{S-KdV1} under
condition that $0<\lambda \leq 1$ and the following condition
\begin{align}\label{indata}
\norm{\phi_1}_{L^2}\leq R,\quad
\norm{\phi_2}_{H^{-3/4}}=\epsilon_0\ll 1.
\end{align}
where $\epsilon_0$ is a absolutely constant will be defined later.
We will prove well-posedness for \eqref{S-KdV1}-\eqref{indata} in
$[0,T]$ for some $T=T(R, \lambda)>0$. By the scaling we obtain local well-posedness
for the original system \eqref{S-KdV}.

In our proof the condition $0<\lambda\leq 1$ in \eqref{S-KdV1} is
crucial. Heuristically, the propagation speed for the first equation
is $\lambda \xi$, and that for the second equation is $\xi^2$. Then
we see that the two waves $u,v$ has a separate speed in high
frequency $|\xi|\ges 1$ uniformly for $0<\lambda\leq 1$. Thus the
resonance and coherence can not be simultaneously large (Also see
\cite{cl07} for the case $\lambda=1$). This is key to control the
coupled wave interactions. Technically, we will prove uniform
estimates for the multiplier associated to the coupled terms for all
$0<\lambda\leq 1$. Our proof for the coupled terms is different from
those in \cite{cl07}, but with basically the same ideas. We will use
the ideas developed by Tao \cite{Taokz}, Ionescu and Kenig
\cite{In-Ke}, and the first-named author \cite{dgGuo}, but in this
paper we need to deal with two different wave forms which has
independent interest.

At the end of this section we introduce some notations. In Section
2, we prove some $L^2$ bilinear estimates which will be used to
prove the bilinear estimates for the coupling terms in Section 3. In
Section 4, we prove Theorem \ref{thmlwp}.

\subsection*{Notations}
Throughout this paper we fix $0<\lambda\leq 1$. We will use $C$ and $c$
to denote constants which are independent of $\lambda$ and not
necessarily the same at each occurrence. For $x, y\in \R$, $x\sim y$
means that there exist $C_1, C_2
> 0$ such that $C_1|x|\leq |y| \leq C_2|x|$. For $f\in \Sch'$ we
denote by $\widehat{f}$ or $\ft (f)$ the Fourier transform of $f$
for both spatial and time variables,
\begin{eqnarray*}
\widehat{f}(\xi, \tau)=\int_{\R^2}e^{-ix \xi}e^{-it \tau}f(x,t)dxdt.
\end{eqnarray*}
We denote  by $\ft_x$ the Fourier transform on spatial variable and
if there is no confusion, we still write $\ft=\ft_x$. Let
$\mathbb{Z}$ and $\mathbb{N}$ be the sets of integers and natural
numbers, respectively. $\Z_+=\N \cup \{0\}$. For $k\in \Z_+$ let
\[{I}_k=\{\xi: |\xi|\in [2^{k-1}, 2^{k+1}]\}, \ k\geq 1; \quad I_0=\{\xi: |\xi|\leq 2\}.\]
Let $\eta_0: \R\rightarrow [0, 1]$ denote an even smooth function
supported in $[-8/5, 8/5]$ and equal to $1$ in $[-5/4, 5/4]$. For
$k\in \Z$ let $\eta_k(\xi)=\eta_0(\xi/2^k)-\eta_0(\xi/2^{k-1})$ if
$k\geq 1$ and $\eta_k(\xi)\equiv 0$ if $k\leq -1$, and let
$\chi_k(\xi)=\eta_0(\xi/2^k)-\eta_0(\xi/2^{k-1})$. For $k\in \Z$ let
$P_k$ denote the operator on $L^2(\R)$ defined by
\[
\widehat{P_ku}(\xi)=\eta_k(\xi)\widehat{u}(\xi).
\]
By a slight abuse of notation we also define the operator $P_k$ on
$L^2(\R\times \R)$ by the formula $\ft(P_ku)(\xi,
\tau)=\eta_k(\xi)\ft (u)(\xi, \tau)$. For $l\in \Z$ let
\[
P_{\leq l}=\sum_{k\leq l}P_k, \quad P_{\geq l}=\sum_{k\geq l}P_k.
\]
Thus we see that $P_{\leq 0}=P_0$.

For $\phi\in \Sch'(\R)$, we denote by
$V(t)\phi=e^{-t\partial_x^3}\phi$ the free solution of linear Airy
equation which is defined as
\[
\ft_x(V(t)\phi)(\xi)=\exp[i\xi^3t]\widehat{\phi}(\xi), \  \forall \
t\in \R,
\]
denote by $U_\lambda(t)\phi=e^{it\lambda\partial_x^2}\phi$ for
$0<\lambda\leq 1$ the free solution of scaled linear Schr\"odinger
equation which is defined as
\[
\ft_x(U_\lambda(t)\phi)(\xi)=\exp[-i\lambda\xi^2t]\widehat{\phi}(\xi),
\ \forall \ t\in \R.
\]
 We define the Lebesgue spaces $L_{t\in I}^qL_x^p$
and $L_x^pL_{t\in I}^q$ by the norms
\begin{equation}
\norm{f}_{L_{t\in I}^qL_x^p}=\normo{\norm{f}_{L_x^p}}_{L_t^q(I)},
\quad \norm{f}_{L_x^pL_{t\in
I}^q}=\normo{\norm{f}_{L_t^q(I)}}_{L_x^p}.
\end{equation}
If $I=\R$ we simply write $L_{t}^qL_x^p$ and $L_x^pL_{t}^q$.

We will make use of the $X^{s,b}$-type space. Generally, let
$h(\xi)$ be a continuous function, and we define
\[\norm{v}_{X_{\tau=h(\xi)}^{s,b}}=\norm{\jb{\tau-h(\xi)}^b\jb{\xi}^s\widehat{v}(\xi,\tau)}_{L^2(\R^2)}\]
where $\jb{\cdot}=(1+|\cdot|^2)^{1/2}$. This type of space was first
systematically studied by Bourgain \cite{b93}. In applications we
usually apply $X^{s,b}$ space for $b$ close to $1/2$. In the case
$b=1/2$ one has a good substitute: Besov-type $X^{s,b}$ space which
was first noted by Tataru \cite{Tataru}. For $k\in \Z_+$ we define
the frequency dyadically localized $X^{s,b}$-type normed spaces
$Y_{\tau=h(\xi)}^k$:
\begin{eqnarray}
Y_{\tau=h(\xi)}^k=\left\{f\in L^2(\R^2):
\begin{array}{l}
f(\xi,\tau) \mbox{ is supported in } I_k\times\R \mbox{ and }\\
\norm{f}_{Y_{\tau=h(\xi)}^k}=\sum_{j=0}^\infty
2^{j/2}\norm{\eta_j(\tau-h(\xi))\cdot f}_{L^2}.
\end{array}
\right\}.
\end{eqnarray}
Then we define the $l^1$-analogue of $X^{s,b}$-type space
$F_{\tau=h(\xi)}^{s}$ by
\begin{align}
\norm{u}_{F_{\tau=h(\xi)}^{s}}^2=\sum_{k \geq
0}2^{2sk}\norm{\eta_k(\xi)\ft(u)}_{Y_{\tau=h(\xi)}^k}^2.
\end{align}
In this paper we will use the space $X_{\tau=-\lambda \xi^2}^{s,b}$
and $F_{\tau=\xi^3}^{s}$. In order to avoid some logarithmic
divergence, we use the following weaker norm for the low frequency
of the KdV equation as in \cite{G},
\begin{align*}
&\norm{u}_{\bar{Y}^0}=\norm{u}_{L_x^2L_t^\infty}.
\end{align*}
For $-3/4\leq s\leq 0$, we define
\begin{align*}
&\bar{F}_{\tau=\xi^3}^s=\{u\in
\Sch'(\R^2):\norm{u}_{\bar{F}_{\tau=\xi^3}^s}^2=\sum_{k \geq
1}2^{2sk}\norm{\eta_k(\xi)\ft(u)}_{Y_{\tau=\xi^3}^k}^2+\norm{P_{\leq
0}(u)}_{\bar{Y}^0}^2<\infty\}.
\end{align*}
For $T\geq 0$, we define the time-localized spaces
$\bar{F}_{\tau=\xi^3}^s(T)$:
\begin{align*}
&\norm{u}_{\bar{F}_{\tau=\xi^3}^s(T)}=\inf_{w\in
\bar{F}_{\tau=\xi^3}^s}\{\norm{P_{\leq 0}u}_{L_x^2L_{|t|\leq
T}^\infty}+\norm{P_{\geq 1}w}_{\bar{F}_{\tau=\xi^3}^s}: \ w(t)=u(t)
\mbox{ on } [-T, T]\}.
\end{align*}
Similarly we define $X_{\tau=h(\xi)}^{s,b}(T)$.

\section{$L^2$ bilinear estimates}

In this section we prove some $L^2$ bilinear estimates which will be
used to prove bilinear estimates for the coupled terms. For
$\xi_1,\xi_2 \in \R$ let
\begin{align}\label{eq:reso}
\Omega_1(\xi_1,\xi_2)=&-\lambda\xi_1^2+\xi_2^3+\lambda(\xi_1+\xi_2)^2,\\
\Omega_2(\xi_1,\xi_2)=&-\lambda\xi_1^2+\lambda\xi_2^2-(\xi_1+\xi_2)^3.
\end{align}
$\Omega_1$ is the resonance function for the coupled term $uv$, and
$\Omega_2$ is the one for $\partial_x(|u|^2)$. For compactly
supported nonnegative functions $f,g,h\in L^2(\R\times \R)$ we
define for $m= 1,2$
\begin{align*}
&J_m(f,g,h)=\int_{\R^4}f(\xi_1,\mu_1)g(\xi_2,\mu_2)h(\xi_1+\xi_2,\mu_1+\mu_2+\Omega_m(\xi_1,\xi_2))d\xi_1d\xi_2d\mu_1d\mu_2.
\end{align*}

For $k,j\in \Z_+$ we define
$$
{D}^\lambda_{k,j}=\{(\xi,\tau):\xi \in I_{k}, \tau+\lambda\xi^2\in
I_{j}\},
$$
and for $k\in \Z, j\in \Z_+$ we define
$$
{B}_{k,j}=\{(\xi,\tau):|\xi| \in [2^{k-1},2^{k+1}], \tau-\xi^3\in
I_{j}\}.
$$

Let $a_1, a_2, a_3\in \R$. It will be convenient to define the
quantities $a_{max}\geq a_{med}\geq a_{min}$ to be the maximum,
median, and minimum of $a_1,a_2,a_3$ respectively. Usually we use
$k_1,k_2,k_3$ and $j_1,j_2,j_3$ to denote integers, $N_i=2^{k_i}$
and $L_i=2^{j_i}$ for $i=1,2,3$ to denote dyadic numbers.

We prove the following lemma.

\begin{lemma}\label{lemsymes}
Assume $k_i \in \Z$, $j_i\in \Z_+$, and $f_{k_i,j_i}\in L^2(\R\times
\R)$ are nonnegative functions supported in
$[2^{k_i-1},2^{k_i+1}]\times I_{j_i}, \, i=1,\ 2,\ 3$. Then

(a) For any $k_1, k_2, k_3\in \Z$ and $j_1,j_2,j_3\in \Z_+$,
$m=1,2$,
\begin{align}\label{eq:lemsymesra}
J_m(f_{k_1,j_1},f_{k_2,j_2},f_{k_3,j_3})&\leq C
2^{j_{min}/2}2^{k_{min}/2} \prod_{i=1}^3\norm{f_{k_i,j_i}}_{L^2}.
\end{align}

(b) If $ k_2\geq 3$, then
\begin{align}\label{eq:lemsymesrb}
J_1(f_{k_1,j_1},f_{k_2,j_2},f_{k_3,j_3})&\leq C
2^{(j_{med}+j_{max})/2}2^{-k_{2}}\prod_{i=1}^3\norm{f_{k_i,j_i}}_{L^2}.
\end{align}

(c) For any $k_1, k_2, k_3\in \Z$ and $j_1,j_2,j_3\in \Z_+$,
\begin{align}\label{eq:lemsymesrc}
J_2(f_{k_1,j_1},f_{k_2,j_2},f_{k_3,j_3})&\leq C\lambda^{-1/2}
2^{(j_{1}+j_{2})/2}2^{- k_{3}/2}
\prod_{i=1}^3\norm{f_{k_i,j_i}}_{L^2}.
\end{align}
\end{lemma}

\begin{proof}
Let $A_{k_i}(\xi)=[\int_\R |f_{k_i,j_i}(\xi,\mu)|^2d\mu]^{1/2}$,
$i=1,2,3$. Using the Cauchy-Schwartz inequality and the support
properties of the functions $f_{k_i,j_i}$, we obtain
\begin{align*}
J_m(f_{k_1,j_1},f_{k_2,j_2},f_{k_3,j_3})\les& 2^{j_{min}/2}\int_{\R^2}A_{k_1}(\xi_1)A_{k_2}(\xi_2)A_{k_3}(\xi_1+\xi_2)d\xi_1d\xi_2\\
\les&
2^{k_{min}/2}2^{j_{min}/2}\prod_{i=1}^3\norm{f_{k_i,j_i}}_{L^2}.
\end{align*}
which is part (a), as desired.

For part (b), in view of the support properties of the functions, it
is easy to see that $J_1(f_{k_1,j_1},f_{k_2,j_2},f_{k_3,j_3})\equiv
0$ unless
\begin{equation}\label{eq:freeq}
|k_{max}-k_{med}|\leq 5, \, 2^{j_{max}}\ges |\Omega_1|.
\end{equation}
We define two sets
$$
A=\{(\xi_1,\xi_2)\in\R^2:
|2\lambda\xi_1-3\xi_2^2|\geq\frac{1}{2}|\xi_2|^2\}
$$
and
$$
B=\{(\xi_1,\xi_2):
|\Omega_1(\xi_1,\xi_2)|\geq\frac{1}{2}|\xi_2|^3\}=\{(\xi_1,\xi_2):
|\xi_2^2+2\lambda\xi_1+\lambda\xi_2|\geq\frac{1}{2}|\xi_2|^2\},
$$
since
$|\Omega_1(\xi_1,\xi_2)|=|\xi_2(\xi_2^2+2\lambda\xi_1+\lambda\xi_2)|$.
We claim that \[[2^{k_1-1},2^{k_1+1}]\times
[2^{k_2-1},2^{k_2+1}]\subset A\cup B.\] Indeed, if
$(\xi_1,\xi_2)\notin A\cup B$, then
$$
|4\xi_2^2+\lambda\xi_2|\leq|3\xi_2^2-2\lambda\xi_1|+|\xi_2^2+2\lambda\xi_1+\lambda\xi_2|\leq
|\xi_2|^2,
$$
which is a contradiction since $|\xi_2|\geq 2$ and $0<\lambda \leq
1$.

For simplicity of notations we set $f_i=f_{k_i,j_2}$, $i=1,2,3$.
Then we get
\begin{align*}
&J_1(f_{1},f_{2},f_{3})\\
\les&\int_{\R^2}\int_Af_{1}(\xi_1,\mu_1)f_{2}(\xi_2,\mu_2)f_{3}(\xi_1+\xi_2,\mu_1+\mu_2+\Omega_m(\xi_1,\xi_2))d\xi_1d\xi_2d\mu_1d\mu_2\\
&+\int_{\R^2}\int_Bf_{1}(\xi_1,\mu_1)f_{2}(\xi_2,\mu_2)f_{3}(\xi_1+\xi_2,\mu_1+\mu_2+\Omega_m(\xi_1,\xi_2))d\xi_1d\xi_2d\mu_1d\mu_2\\
=&I+II.
\end{align*}
We consider first the contribution of the term $I$. Since for
$(\xi_1,\xi_2)\in A$, then
\[|\partial_{\xi_1}\Omega_1-\partial_{\xi_2}\Omega_1|=|2\lambda\xi_1-3\xi_2^2|\geq\frac{1}{2}|\xi_2|^2.\]
We will prove that if $g_i:\R\rightarrow \R_+$ are $L^2$ functions
supported in $[2^{k_i-1},2^{k_i+1}]$, $i=1,2$, and $g:
\R^2\rightarrow \R_+$ is an $L^2$ function supported in
$[2^{k_3-1},2^{k_3+1}]\times I_{j_{max}}$, then
\begin{eqnarray}\label{eq:lemsymesb}
\int_{A}g_1(\xi_1)g_2(\xi_2)g(\xi_1+\xi_2,\Omega_1(\xi_1,\xi_2))d\xi_1d\xi_2\les
2^{-k_2}\norm{g_1}_{L^2}\norm{g_2}_{L^2}\norm{g}_{L^2}.
\end{eqnarray}
This suffices for \eqref{eq:lemsymesrb} by Cauchy-Schwartz
inequality. To prove \eqref{eq:lemsymesb} we get
\begin{align*}
&\int_{A}g_1(\xi_1)g_2(\xi_2)g(\xi_1+\xi_2,\Omega_1(\xi_1,\xi_2))d\xi_1d\xi_2\\
&\les\norm{g_1}_{L^2}\norm{g_2}_{L^2}\norm{1_{A}(\xi_1,\xi_2)g(\xi_1+\xi_2,\Omega_1(\xi_1,\xi_2))}_{L^2_{\xi_1\xi_2}}\\
&\les2^{-k_2}\norm{g_1}_{L^2}\norm{g_2}_{L^2}\norm{g}_{L^2},
\end{align*}
where in the last inequality we used the change of the variables
$x=\xi_1+\xi_2, y=\Omega_1(\xi_1,\xi_2)$, since the Jacobi is
$|\partial_{\xi_1}\Omega_1-\partial_{\xi_2}\Omega_1|\geq\frac{1}{2}|\xi_2|^2\ges
2^{2k_2}$ in $A$.

Now we consider the term $II$. When $(\xi_1,\xi_2)\in B$, we have
that $|\Omega_1|\geq \frac{1}{2}|\xi_2|^3$, so from the support
properties we get $j_{max}\geq 3k_2-20$. Then from (a) we have in
this case
\[
II\les2^{j_{min}/2}2^{k_{2}/2}
\prod_{i=1}^3\norm{f_{k_i,j_i}}_{L^2}\les2^{j_{min}/2}2^{j_{max}/2}2^{-k_2}
\prod_{i=1}^3\norm{f_{k_i,j_i}}_{L^2},
\]
which completes the proof of (b).

Now we prove part (c). We will prove that if $g_i:\R\rightarrow
\R_+$ are $L^2$ functions supported in $[2^{k_i-1},2^{k_i+1}]$,
$i=1,2$, and $g: \R^2\rightarrow \R_+$ is an $L^2$ function
supported in $[2^{k_3-1},2^{k_3+1}]\times I_{j_{max}}$, then
\begin{eqnarray}\label{eq:lemsymesb2}
\int_{\R^2}g_1(\xi_1)g_2(\xi_2)g(\xi_1+\xi_2,\Omega_2(\xi_1,\xi_2))d\xi_1d\xi_2\les\lambda^{-1/2}2^{-k_3/2}\norm{g_1}_{L^2}\norm{g_2}_{L^2}\norm{g}_{L^2}.
\end{eqnarray}
This suffices for \eqref{eq:lemsymesrc} by Cauchy-Schwartz
inequality. To prove \eqref{eq:lemsymesb2}, we notice that
\[|\partial_{\xi_1}\Omega_2-\partial_{\xi_2}\Omega_2|=2\lambda|\xi_1+\xi_2|\thicksim
2\lambda2^{k_3},\] thus by change of variable $\mu_1=\xi_1+\xi_2$,
$\mu_2=\Omega_2(\xi_1,\xi_2)$ we get
\[\|g(\xi_1+\xi_2,\Omega_2(\xi_1,\xi_2))\|_{L^2_{\xi_1\xi_2}}\les
\lambda^{-1/2}2^{-k_3/2}\norm{g}_2,\]which is sufficient for
\eqref{eq:lemsymesb2} by Cauchy-Schwartz inequality.
\end{proof}

\begin{remark}\label{rm:symest}
It is easy to see from the proof that Part (a) and Part (b) of Lemma
\ref{lemsymes} also hold if we assume instead $f_{k_i,j_i}$ is
supported in $I_{k_i}\times I_{j_i}$ for $k_1,k_2,k_3\in \Z_+$. Part
(c) of Lemma \ref{lemsymes} also holds if we assume instead
$f_{k_1,j_1}$, $f_{k_2,j_2}$ are supported in $I_{k_1}\times
I_{j_1}$, $I_{k_2}\times I_{j_2}$ respectively, for $k_1,k_2\in
\Z_+$.
\end{remark}

We restate now Lemma \ref{lemsymes} in a form that is suitable for
the bilinear estimates in the next section.

\begin{corollary}\label{cor42}

(a) Let $k_1,k_2,k_3,j_1,j_2,j_3\in \Z_+$. Assume
$f_{k_1,j_1},f_{k_2,j_2}\in L^2(\R\times \R)$ are nonnegative
functions that are supported in $\{(\xi,\tau):\xi \in I_{k_1},
\tau+\lambda\xi^2\in I_{j_1}\}$ and $\{(\xi,\tau):\xi \in I_{k_2},
\tau-\xi^3\in I_{j_2}\}$ respectively, then
\begin{eqnarray}\label{aa}
\norm{1_{{D}^\lambda_{k_3,j_3}}(\xi,\tau)(f_{k_1,j_1}*f_{k_2,j_2})}_{L^2}\les
2^{k_{min}/2}2^{j_{min}/2}\prod_{i=1}^2\norm{f_{k_i,j_i}}_{L^2}.
\end{eqnarray}
Furthermore, if $k_2>2$, then we have
\begin{eqnarray}\label{k2}
\norm{1_{{D}^\lambda_{k_3,j_3}}(\xi,\tau)(f_{k_1,j_1}*f_{k_2,j_2})}_{L^2}\les
2^{(j_{med}+j_{max})/2}2^{-
k_{2}}\prod_{i=1}^2\norm{f_{k_i,j_i}}_{L^2}.
\end{eqnarray}

(b) Let $k_1,k_2, j_1,j_2,j_3\in \Z_+$ and $k_3\in \Z$. Assume
$f_{k_1,j_1},f_{k_2,j_2}\in L^2(\R\times \R)$ are nonnegative
functions that are supported in $\{(\xi,\tau):\xi \in {I}_{k_1},
\tau+\lambda\xi^2\in I_{j_1}\}$ and $\{(\xi,\tau):\xi \in I_{k_2},
\tau-\lambda\xi^2\in I_{j_2}\}$ respectively, then
\begin{eqnarray}\label{bb}
\norm{1_{{B}_{k_3,j_3}}(\xi,\tau)(f_{k_1,j_1}*f_{k_2,j_2})}_{L^2}\les
2^{k_{min}/2}2^{j_{min}/2}\prod_{i=1}^2\norm{f_{k_i,j_i}}_{L^2},
\end{eqnarray}
and also
\begin{eqnarray}\label{b2}
\norm{1_{{B}_{k_3,j_3}}(\xi,\tau)(f_{k_1,j_1}*f_{k_2,j_2})}_{L^2}\les
\lambda^{-1/2}
2^{j_{1}/2}2^{j_{2}/2}2^{-k_{3}/2}\prod_{i=1}^2\norm{f_{k_i,j_i}}_{L^2}.
\end{eqnarray}
\end{corollary}

\begin{proof} We just prove (a), the proof for (b) is similar.
Clearly, we have
\begin{align*}
&\norm{1_{{D}^\lambda_{k_3,j_3}}(\xi,\tau)(f_{k_1,j_1}*f_{k_2,j_2})(\xi,\tau)}_{L^2}=\sup_{\norm{f}_{L^2}=1}\aabs{\int_{D^\lambda_{k_3,j_3}}
f\cdot f_{k_1,j_1}*f_{k_2,j_2} d\xi d\tau}.
\end{align*}
We denote $f_{k_3,j_3}=1_{D^\lambda_{k_3,j_3}}(\xi,\tau)\cdot
f(\xi,\tau)$. Define
$f_{k_1,j_1}^\sharp(\xi,\mu)=f_{k_1,j_1}(\xi,\mu-\lambda\xi^2)$,
$f_{k_2,j_2}^\sharp(\xi,\mu)=f_{k_2,j_2}(\xi,\mu+\xi^3)$,
$f_{k_3,j_3}^\sharp(\xi,\mu)=f_{k_3,j_3}(\xi,\mu-\lambda\xi^2)$.
Then for $i=1,2,3$ the functions $f_{k_i,j_i}^\sharp$ are supported
in $I_{k_i}\times I_{j_i}$  and
$\norm{f_{k_i,j_i}^\sharp}_{L^2}=\norm{f_{k_i,j_i}}_{L^2}$. Using
simple changes of variables, we get that
\[\int_{D^\lambda_{k_3,j_3}} f\cdot
f_{k_1,j_1}*f_{k_2,j_2} d\xi d\tau =
J_1(f_{k_1,j_1}^\sharp,f_{k_2,j_2}^\sharp,f_{k_3,j_3}^\sharp).\]
Then (a) follows from Lemma \ref{lemsymes} (a), (b) and Remark
\ref{rm:symest}.
\end{proof}

\section{Bilinear estimates for the coupling terms}
This section is devoted to prove the bilinear estimates for the
coupling terms in the $F^s$-type space. We first recall an abstract
extension lemma. For a continuous function $h(\xi)$, define group
$W_h(t)$ by
\[W_h(t)f=\ft_x e^{ith(\xi)}\ft_x{f}.\]
Then $W_h(t)f$ is the solution to the following equation
\[\partial_t u-ih(-i\partial_x)u=0,\quad u(x,0)=f(x).\]
\begin{lemma}[Extension lemma]\label{lemmaextension}
Assume $h$ is a continuous function. Let $Z$ be any space-time
Banach space which obeys the time modulation estimate
\begin{eqnarray}
\norm{g(t)F(t,x)}_{Z}\leq \norm{g}_{L_t^\infty} \norm{F(t,x)}_Z
\end{eqnarray}
for any $F\in Z$ and $g\in L_t^\infty$. Moreover, if for all
$u_{0}\in L_x^2$
\[\norm{W_h(t)u_{0}}_{Z}\les \norm{u_{0}}_{L_x^2}.\]
Then one also has the estimate that for all $k\in \Z_+$ and $u\in
{F}^0$
\[\norm{P_{k}(u)}_{Z}\les \norm{\wh{P_{k}(u)}}_{{Y}_{\tau=h(\xi)}^{k}}.\]
\end{lemma}
\begin{proof}
We refer the readers to Lemma 3.2 in \cite{G}.
\end{proof}

\begin{proposition}[$X_{\tau=-\lambda \xi^2}^k$ embedding]\label{propXkembedding}
Let $k\in \Z_+$, $j\in \N$. Assume $u\in \Sch(\R^2)$, then we have
\begin{align*}
\norm{P_k(u)}_{L_t^qL_x^r}\les&\lambda^{-1/q}\norm{P_k(u)}_{{Y}_{\tau=-\lambda\xi^2}^{k}},\label{str}\\
\norm{P_j(u)}_{L_x^\infty L_{t}^2}\les&
\lambda^{-1/2}2^{-j/2}\norm{\ft[P_j(u)]}_{{Y}_{\tau=-\lambda\xi^2}^{k}},
\end{align*}
where $(q,r)$ satisfies $2\leq q,r \leq \infty$ and 2/q=1/2-1/r. As
a consequence, we get from the definition that for $u\in
F_{\tau=-\lambda\xi^2}^{s}$
\[\norm{u}_{L_t^\infty H^s}\les \norm{u}_{F_{\tau=-\lambda\xi^2}^{s}}.\]
\end{proposition}
\begin{proof}
From Lemma \ref{lemmaextension}, it suffices to prove
\begin{align*}
\norm{U_\lambda(t)f}_{L_t^qL_x^r}\les&\lambda^{-1/q}\norm{f}_{L_x^2},\label{str}\\
\norm{U_\lambda(t)f}_{L_x^\infty
L_{t}^2}\les&\lambda^{-1/2}\norm{f}_{\dot{H}^{-1/2}_x},
\end{align*}
which are well-known, for example see \cite{KeelTao} and
\cite{KPVIUMJ91}.
\end{proof}

\begin{proposition}[$X_{\tau=\xi^3}^k$ embedding]\label{propYkembedding}
Let $k\in \Z_+$, $j\in \N$. Assume $u\in F^{0}$, then we have
\begin{eqnarray}
&&\norm{P_k(u)}_{L_t^qL_x^r}\les \norm{P_k(u)}_{Y_{\tau=\xi^3}^k},\\
&&\norm{P_k(u)}_{L_x^2L_{t\in I}^\infty}\les 2^{3k/4}
\norm{\ft[P_k(u)]}_{Y_{\tau=\xi^3}^k},\\
&&\norm{P_j(u)}_{L_x^\infty L_{t}^2}\les
2^{-j}\norm{\ft[P_j(u)]}_{Y_{\tau=\xi^3}^j},
\end{eqnarray}
where $(q,r)$ satisfies $2\leq q,r \leq \infty$ and 3/q=1/2-1/r. As
a consequence, we get from the definition that for $u\in
\bar{F}_{\tau=\xi^3}^s$, then
\[\norm{u}_{L_t^\infty H^s}\les \norm{u}_{\bar{F}_{\tau=\xi^3}^s}.\]
\end{proposition}
\begin{proof}
We refer the readers to proposition 3.3 in \cite{G}.
\end{proof}

The main result of this section is the following lemma.

\begin{lemma}\label{biline}
(a) If $s_1\geq 0$, $s_2\in (-1,-1/2)$, $s_1-s_2<1$, $0<\theta\ll
1$, and $u\in X^{s_1,1/2+\theta}_{\tau=-\lambda\xi^2}$, $v\in
\bar{F}_{\tau=\xi^3}^{s_2}$ then
\begin{align}\label{eq:lemcouplea}
\norm{\psi(t)uv}_{X^{s_1,-1/2+2\theta}_{\tau=-\lambda\xi^2}}\les&\lambda^{-1/2}
\norm{u}_{X^{s_1,1/2+\theta}_{\tau=-\lambda\xi^2}}\norm{v}_{\bar{F}_{\tau=\xi^3}^{s_2}}.
\end{align}

(b)If $s_1\geq 0$, $s_2-s_1<-1/2$, and $u, w\in
X^{s_1,1/2+\theta}_{\tau=-\lambda\xi^2}$ then
\begin{align}\label{eq:lemcoupleb}
\norm{\psi(t)\partial_x(u\bar{w})}_{{X}_{\tau=\xi^3}^{s_2,-1/2+2\theta}}\les&
\lambda^{-1/2}\norm{u}_{X^{s_1,1/2+\theta}_{\tau=-\lambda\xi^2}}\norm{w}_{X^{s_1,1/2+\theta}_{\tau=-\lambda\xi^2}}.
\end{align}
\end{lemma}

\begin{proof}
We first prove (a). From the definition, we get that the left-hand
side of \eqref{eq:lemcouplea} is dominated by
\begin{align}\label{pa1}
\Big\|2^{s_1 k_3}\sum_{j_3=0}^{\infty}2^{-1/2j_3+2\theta
j_3}\|1_{D^\lambda_{k_3,j_3}}\cdot\widehat{\psi(t)uv}\|_{L_{\xi,\tau}^{2}}\Big\|_{l^{2}_{k_3}}.
\end{align}
Now we begin to estimate
\begin{align}\label{pa1}
\sum_{j_3}^{\infty}2^{-1/2j_3+2\theta
j_3}\|1_{D^\lambda_{k_3,j_3}}\cdot\widehat{\psi(t)uv}\|_{L_{\xi,\tau}^{2}}.
\end{align}
Decomposing $u, v$, for $k_1,j_1\in \Z_+$ set
\begin{align*}
f_{k_1,j_1}(\xi_1,\tau_1)&=\eta_{k_1}(\xi_1)\eta_{j_1}(\tau_1+\lambda\xi_1^2)\widehat{\psi(t/2)u}(\xi_1,\tau_1),\\
f_{k_2,j_2}(\xi_2,\tau_2)&=\eta_{k_2}(\xi_2)\eta_{j_2}(\tau_2-\xi_2^3)\widehat{\psi(t/2)v}(\xi_2,\tau_2),
\end{align*}
then we get for fixed $k_3$,
\begin{eqnarray}\label{pa2}
\eqref{pa1}\les\sum_{(k_1,k_2,k_3)\in K } \sum_{j_3\geq
0}2^{-j_3/2+2\theta j_3} \sum_{j_1, j_2\geq
0}\|1_{D_{k_3,j_3}^\lambda}\cdot
f_{k_1,j_1}*f_{k_2,j_2}\|_{L^2_{\xi_3\tau_3}}.
\end{eqnarray}
where $K=\{(k_1,k_2,k_3)\in \Z_+^3: |k_{med}-k_{max}|<5\}$, since it
is easy to see that $1_{D_{k_3,j_3}^\lambda}\cdot
f_{k_1,j_1}*f_{k_2,j_2}\equiv 0$ unless
\[|k_{med}-k_{max}|<5.\]
We may also assume that $k_{max} \geq 20$, since for $k_{max} \leq
20$ we can get from Plancherel's equality that
\begin{align*}
\eqref{pa1}\les&\sum_{j_3\geq
0}2^{-j_3/2+2\theta j_3}\|1_{D_{k_3,j_3}^\lambda}\cdot\wh{P_{\leq 20}u}*\wh{\psi(t)P_{\leq 20}v}\|_{L^2_{\xi_3\tau_3}}\\
\les&\|{P_{\leq 20}u}{\psi(t)P_{\leq 20}v}\|_{L^2_{x,t}}\\
\les& \norm{P_{\leq 20}u}_{L_t^\infty L_x^2}\norm{\psi(t)P_{\leq
20}v}_{ L_x^2L_t^\infty}
\end{align*}
which suffices to give the bound for this case by Proposition
\ref{propXkembedding}, \ref{propYkembedding}.

Now we assume $k_{max}\geq 20$ and prove \eqref{pa2}. First we
assume that $k_2\geq 2$. Clearly we may also assume that
$j_{max}\leq 10 k_2$, otherwise, we can apply \eqref{aa}, then we
have a $2^{-5k_2}$ to spare. After these assumptions, we can make
use of \eqref{k2} to bound \eqref{pa2} by
\begin{eqnarray}\label{pa3}
\sum_{(k_1,k_2,k_3)\in K, k_2\geq 2} \sum_{j_3\geq
0}2^{-j_3/2+2\theta j_2} \sum_{j_1, j_2\geq
0}2^{j_{max}/2}2^{j_{med}/2}2^{-k_{2}}\prod_{i=1}^2\norm{f_{k_i,j_i}}_{L^2}
\end{eqnarray}
when $j_3=j_{min}$, then \eqref{pa3} is bounded by
\begin{align}\label{pa31}
&\sum_{(k_1,k_2,k_3)\in K, k_2\geq 2} \sum_{j_1, j_2\geq
0}2^{j_{1}/2}2^{j_{2}/2}2^{-k_{2}}\prod_{i=1}^2\norm{f_{k_i,j_i}}_{L^2}\nonumber\\
&\les\sum_{(k_1,k_2,k_3)\in K, k_2\geq 2}
2^{-k_{2}}\norm{P_{k_1}u}_{Y_{\tau=-\lambda\xi^2}^{k_1}}\norm{P_{k_2}v}_{Y_{\tau=\xi^3}^{k_2}},
\end{align}
when $j_1=j_{min}$, then \eqref{pa3} is bounded by
\begin{align}\label{pa32}
&\sum_{(k_1,k_2,k_3)\in K, k_2\geq 2}\sum_{j_3\geq 0}2^{2\theta j_3}
\sum_{j_1, j_2\geq
0}2^{j_{2}/2}2^{-k_{2}}\prod_{i=1}^2\norm{f_{k_i,j_i}}_{L^2}\nonumber\\
&\les\sum_{(k_1,k_2,k_3)\in K, k_2\geq 2}2^{20\theta k_2} \sum_{j_1,
j_2\geq
0}2^{j_{1}/2}2^{j_{2}/2}2^{-k_{2}}\prod_{i=1}^2\norm{f_{k_i,j_i}}_{L^2}\nonumber\\
&\les\sum_{(k_1,k_2,k_3)\in K, k_2\geq 2} 2^{-k_{2}+20\theta
k_2}\norm{P_{k_1}u}_{Y_{\tau=-\lambda\xi^2}^{k_1}}\norm{P_{k_2}v}_{Y_{\tau=\xi^3}^{k_2}},
\end{align}
where we use $j_{max}\leq 10k_2$ in the last step. The case
$j_2=j_{min}$ is the same as $j_1=j_{min}$, so we omit the details.

Dividing the summation on $k_1,k_2$ in the right-hand side of
\eqref{pa3} into several parts, we get from
\eqref{pa31},\eqref{pa32}
\begin{align*}
&\eqref{pa3}\les \sum_{i=1}^5 \sum_{{A_i(k_3)}}2^{-k_{2}+20\theta
k_2}\norm{P_{k_1}u}_{Y_{\tau=-\lambda\xi^2}^{k_1}}\norm{P_{k_2}v}_{Y_{\tau=\xi^3}^{k_2}}
\end{align*}
where we denote
\begin{align*}
A_1(k_3)=&\{|k_2-k_3|\leq 5, k_1\leq k_2-10, \mbox{ and } k_2\geq
30\};\\
A_2(k_3)=&\{|k_1-k_3|\leq 5, 2\leq k_2\leq k_1-10, \mbox{ and }
k_1\geq
30\};\\
A_3(k_3)=&\{|k_1-k_2|\leq 5,  k_3\leq k_2-10, \mbox{ and } k_1 \geq
30\};\\
A_4(k_3)=&\{|k_2-k_3|\leq 10, |k_1-k_2|\leq 10, \mbox{ and } k_2\geq
30 \};\\
A_5(k_3)=&\{k_1,k_2,k_3\leq 200\}.
\end{align*}
So we can bounded \eqref{pa1} by
\begin{align}\label{pa11}
\sum_{i=1}^5 \Big\|2^{s_1 k_3}\sum_{{A_i(k_3)}}2^{-k_{2}+20\theta
k_2}\norm{P_{k_1}u}_{Y_{\tau=-\lambda\xi^2}^{k_1}}\norm{P_{k_2}v}_{Y_{\tau=\xi^3}^{k_2}}\Big\|_{l^{2}_{k_3}}.
\end{align}

Now we begin to estimate \eqref{pa11} case by case. For case $A_1$,
by Cauchy-Schwartz's inequality we have
\begin{align}\label{pa1a1}
&\Big\|2^{s_1 k_3}\sum_{{A_1(k_3)}}2^{-k_{2}+20\theta
k_2}\norm{P_{k_1}u}_{Y_{\tau=-\lambda\xi^2}^{k_1}}\norm{P_{k_2}v}_{Y_{\tau=\xi^3}^{k_2}}\Big\|_{l^{2}_{k_3}}\nonumber\\
=&\Big\|2^{s_1 k_3}\sum_{|k_2-k_3|\leq 5, k_2\geq 20}\sum_{k_1\leq
k_2-10}2^{-k_{2}+20\theta
k_2}\norm{P_{k_1}u}_{Y_{\tau=-\lambda\xi^2}^{k_1}}\norm{P_{k_2}v}_{Y_{\tau=\xi^3}^{k_2}}\Big\|_{l^{2}_{k_3}}\nonumber\\
\les&\Big\|2^{s_1 k_3}\sum_{|k_2-k_3|\leq 5, k_2\geq
20}2^{-k_{2}+21\theta
k_2}\norm{P_{k_2}v}_{Y_{\tau=\xi^3}^{k_2}}\Big\|_{l^{2}_{k_3}}\Big(\sum_{k_1\geq
0
}2^{2s_1k_{1}}\norm{P_{k_1}u}_{Y_{\tau=-\lambda\xi^2}^{k_1}}\Big)^{1/2}.
\end{align}
since $s_1-s_2<1$, $s_1\geq0$, $s_2>-1$ and $0<\theta\ll 1$, we can
bound \eqref{pa1a1} by
$\|v\|_{\bar{F}_{\tau=\xi^3}^{s_2}}\|u\|_{X^{s_1,1/2+\theta}_{\tau=\lambda\xi^2}}$.
The proof for $A_2-A_5$ are similar, we omit here.

Now we assume that $k_2\leq 2$, by the assumption $k_{max}\geq 20$,
 we have $|k_1-k_3|\leq 5$, and $k_1, k_2 \geq 10$. For simplicity of notations we assume
 $k_1=k_3=k$, then the left hand
side of \eqref{pa1} can be estimated as
\begin{eqnarray}
\sum_{k\geq 0} 2^{s_1k} \|\widehat{P_k(u)}*\widehat{\psi(t)P_{\leq
2}(v)}\|_{L^2}\les \sum_{k\geq 0} 2^{s_1k} \|P_k(u)\|_{L_x^\infty
L_t^2}\|\psi(t)P_{\leq 2}(v)\|_{L_x^2 L_t^\infty}.
\end{eqnarray}
This is enough in view of the definition, Proposition
\ref{propXkembedding} and Proposition \ref{propYkembedding}.

\medskip

Now we begin to prove (b), form the definition, we only need to
prove
\begin{eqnarray}\label{pb}
\norm{\partial_x(u\bar{w})}_{{X}_{\tau=\xi^3}^{s_2,-1/2+2\theta}}\les
\lambda^{-1/2}
\norm{u}_{F^{s_1}_{\tau=-\lambda\xi^2}}\norm{w}_{F^{s_1}_{\tau=-\lambda\xi^2}}.
\end{eqnarray}
By the definition and $\theta \ll 1$, we have
\begin{align}\label{pb1}
\norm{\partial_x(u\bar{w})}_{{X}_{\tau=\xi^3}^{s_2,-1/2+2\theta}}\les&\Big(\sum_{k
\geq 1}2^{2s_2k}\Big(\sum_{j\geq
0}2^{-j/4}\norm{1_{B_{k,j}}(\xi,\tau)i\xi\ft(u\bar{w})}_{L^2_{\xi,\tau}}\Big)^{2}\Big)^{1/2}\nonumber\\
&+\Big(\sum_{k \leq 0}\Big(\sum_{j\geq
0}2^{-j/4}\norm{1_{B_{k,j}}(\xi,\tau)i\xi\ft(u\bar{w})}_{L^2_{\xi,\tau}}\Big)^{2}\Big)^{1/2}\nonumber\\
:=&I+II.
\end{align}
Now let
$$f_{k_1,j_1}(\xi_1,\tau_1)=\eta_{k_1}(\xi_1)\eta_{j_1}(\tau_1+\lambda\xi_1^2)\hat{u}(\xi_1,\tau_1),$$
$$\overline{g_{k_2,j_2}(\xi_2,\tau_2)}=\eta_{k_1}(\xi_2)\eta_{j_1}(\tau_2-\lambda\xi_2^2)\widetilde{\hat{w}}(\xi_2,\tau_2),$$
then we have
\[
\|u\|_{F_{\tau=-\lambda\xi^2}^{s}} = \Big\|2^{sk_1}\sum_{j_1\geq
0}2^{j_1/2}\|f_{k_1,j_1}\|_{L^2}\Big\|_{l^2_{k_1}},\]
\[
\|w\|_{F_{\tau=-\lambda\xi^2}^{s}}= \Big\|2^{sk_2}\sum_{j_2\geq
0}2^{j_2/2}\|g_{k_2,j_2}\|_{L^2}\Big\|_{l^2_{k_2}}.
\]
Then from the definition, $I$ and $II$ in \eqref{pb1} can be bounded
as following
\begin{eqnarray}\label{pb2}
I\les \Big\|2^{(s_2+1)k_3}\sum_{(k_1,k_2,k_3)\in K} \sum_{j_3\geq
0}2^{-j_3/4} \sum_{j_1, j_2\geq 0}\|1_{B_{k_3,j_3}}
f_{k_1,j_1}*g_{k_2,j_2}\|_{L^2_{\xi_3\tau_3}}\Big\|_{l_{k_3\geq
1}^2}
\end{eqnarray}
and
\begin{eqnarray}\label{pb3}
II\les\Big\|2^{k_3}\sum_{(k_1,k_2,k_3)\in K} \sum_{j_3\geq
0}2^{-j_3/4} \sum_{j_1, j_2\geq 0}\|1_{B_{k_3,j_3}}
f_{k_1,j_1}*g_{k_2,j_2}\|_{L^2_{\xi_3\tau_3}}\Big\|_{l^2_{ k_3\leq
0}}.
\end{eqnarray}
We first estimate term $I$ in \eqref{pb2}. By \eqref{b2} we have
\begin{align}\label{pb5}
&\sum_{(k_1,k_2,k_3)\in K}\sum_{j_3\geq 0}2^{-j_3/4} \sum_{j_1,
j_2\geq 0}\|1_{B_{k_3,j_3}}
f_{k_1,j_1}*g_{k_2,j_2}\|_{L^2_{\xi_3\tau_3}}\nonumber\\
\les& \sum_{(k_1,k_2,k_3)\in K, k_3\geq 1} 2^{-k_3/2}\sum_{j_1,
j_2\geq 0}\lambda^{-1/2}
2^{j_{1}/2}2^{j_{2}/2}\norm{f_{k_i,j_i}}_{L^2}\norm{g_{k_i,j_i}}_{L^2}\nonumber\\
\les&\lambda^{-1/2}\sum_{(k_1,k_2,k_3)\in K, k_3\geq 1}2^{-k_3/2}
\norm{\wh{P_{k_1}u}}_{Y_{\tau=-\lambda\xi^2}^{k_1}}\norm{\wh{P_{k_2}w}}_{Y_{\tau=-\lambda\xi^2}^{k_2}}.
\end{align}

Dividing the summation on $k_1,k_2$ in the right-hand side of
\eqref{pb5} into several parts,
\begin{align*}
&\eqref{pb5}\les \lambda^{-1/2}\sum_{i=1}^5
\sum_{{B_i(k_3)}}\norm{P_{k_1}u}_{F_{\tau=\lambda\xi^2}^{k_1}}\norm{P_{k_2}w}_{F_{\tau=\lambda\xi^2}^{k_2}}
\end{align*}
where we denote
\begin{align*}
B_1(k_3)=&\{|k_2-k_3|\leq 5, k_1\leq k_2-10, \mbox{ and } k_2\geq
30\};\\
B_2(k_3)=&\{|k_1-k_3|\leq 5,  k_2\leq k_1-10, \mbox{ and } k_1\geq
30\};\\
B_3(k_3)=&\{|k_1-k_2|\leq 5,  1\leq k_3\leq k_2-10, \mbox{ and } k_1
\geq
30\};\\
B_4(k_3)=&\{|k_2-k_3|\leq 10, |k_1-k_2|\leq 10, \mbox{ and } k_2\geq
30 \};\\
B_5(k_3)=&\{k_1,k_2,k_3\leq 200, k_3\geq 1\}.
\end{align*}
So we can bounded \eqref{pb2} by
\begin{align}\label{pb11}
I\les \lambda^{-1/2}\sum_{i=1}^5
\Big\|2^{(s_2+1/2)k_3}\sum_{{B_i(k_3)}}\norm{\wh{P_{k_1}u}}_{Y_{\tau=-\lambda\xi^2}^{k_1}}\norm{\wh{P_{k_2}w}}_{Y_{\tau=-\lambda\xi^2}^{k_2}}\Big\|_{l^{2}_{k_3}}.
\end{align}

Now we begin to estimate \eqref{pb11} case by case. For case $B_1$,
by Cauchy-Schwartz's inequality we have
\begin{align}\label{pb1a1}
&\lambda^{-1/2}\Big\|2^{(s_2+1/2)k_3}\sum_{{B_1(k_3)}}\norm{\wh{P_{k_1}u}}_{Y_{\tau=-\lambda\xi^2}^{k_1}}\norm{\wh{P_{k_2}w}}_{Y_{\tau=-\lambda\xi^2}^{k_2}}\Big\|_{l^{2}_{k_3}}\nonumber\\
=&\lambda^{-1/2}\Big\|2^{(s_2+1/2)k_3}\sum_{|k_2-k_3|\leq 5, k_2\geq
20}\norm{\wh{P_{k_2}w}}_{Y_{\tau=-\lambda\xi^2}^{k_2}}\sum_{k_1\leq k_2-10}\norm{\wh{P_{k_1}u}}_{Y_{\tau=-\lambda\xi^2}^{k_1}}\Big\|_{l^{2}_{k_3}}\nonumber\\
\les&\lambda^{-1/2}\Big\|2^{(s_2+1/2+\theta)k_3}\sum_{|k_2-k_3|\leq
5, k_2\geq
20}\norm{\wh{P_{k_2}w}}_{Y_{\tau=-\lambda\xi^2}^{k_2}}\Big\|_{l^{2}_{k_3}}\|u\|_{F_{\tau=-\lambda\xi^2}^{s_1}}.
\end{align}
since $s_2-s_1<-1/2$, $s_1\geq0$ and $0<\theta\ll 1$, so we can
bounded \eqref{pb1a1} by
$\|u\|_{X^{s_1,1/2+\theta}_{\tau=-\lambda\xi^2}}\|w\|_{X^{s_1,1/2+\theta}_{\tau=-\lambda\xi^2}}$.
The proof for $B_2-B_5$ are similar, we omit here.

Use the same argument as above, under the restriction
$s_2-s_1<-1/2$, $s_1\geq0$ and $0<\theta\ll 1$, we can bound the
part $II$ in \eqref{pb3} by
$\|u\|_{X^{s_1,1/2+\theta}_{\tau=-\lambda\xi^2}}\|w\|_{X^{s_1,1/2+\theta}_{\tau=-\lambda\xi^2}}$.
Thus we finish the proof.
\end{proof}

\section{Proof of the main theorem}

In this section we prove Theorem \ref{thmlwp}. We first recall some
linear estimates in $X_{\tau=h(\xi)}^{s,b}$ and
$\bar{F}_{\tau=h(\xi)}^s$. Let $W_h(t)f=\ft_x^{-1}e^{ith(\xi)}\ft_x
f$. The following lemma has been proved by Kenig, Ponce and Vega in
\cite{kpvd93}, and then improved by Ginibre, Y. Tsutsumi and Velo in
\cite{gtv97}.

\begin{lemma}\label{line}
Let $s\in \R$, $-1/2<b'\leq 0\leq b\leq b'+1$ and $T\in [0,1]$. Then
for $u_0\in H^s$ and $F\in X^{s,b}_{\tau=h(\xi)}$ we have
\begin{align*}
\|\psi_1(t)W_h(t)u_0\|_{X^{s,b}_{\tau=h(\xi)}}\les&
\|u_0\|_{H^s},\\
\normo{\psi_T(t)\int_{0}^{t}W_h(t-t')F(t',\cdot)dt'}_{X^{s,b}_{\tau=h(\xi)}}\les&
T^{1-b+b'}\|F\|_{X^{s,b'}_{\tau=h(\xi)}}.
\end{align*}
\end{lemma}

Next we prove the linear estimates in $\bar{F}^s$. The proof was
essentially given in \cite{In-Ke} for the Benjamin-Ono equation.
\begin{lemma}
(a) Assume $s\in \R$ and  $\phi \in H^s$. Then there exists $C>0$
such that
\begin{eqnarray}
\norm{\psi(t)W_h(t)\phi}_{{F}_{\tau=h(\xi)}^s}\leq
C\norm{\phi}_{H^{s}}.
\end{eqnarray}

(b) Assume $s\in \R, k\in \Z_+$ and $u$ satisfies
$(i+\tau-\xi^3)^{-1}\ft(u)\in Y_{\tau=h(\xi)}^k$. Then there exists
$C>0$ such that
\begin{eqnarray}
\normo{\ft\left[\psi(t)\int_0^t
W_h(t-s)(u(s))ds\right]}_{Y_{\tau=h(\xi)}^k}\leq
C\norm{(i+\tau-\xi^3)^{-1}\ft(u)}_{Y_{\tau=h(\xi)}^k}.
\end{eqnarray}
\end{lemma}

Now we give the estimates for the cubic nonlinear term.
\begin{lemma}\label{nonS}
Let $u,u'\in X^{s,b}_{\tau=-\lambda\xi^2}$ with $1/2<b<1$ and $s\geq
0$. Then for all $a\geq 0$ we have that
\begin{align*}
\||u|^2u\|_{X^{s,-a}_{\tau=-\lambda\xi^2}}\les& \lambda
^{-1/2}\|u\|^3_{X^{s,b}_{\tau=-\lambda\xi^2}},\\
\||u|^2u-|u'|^2u'\|_{X^{s,-a}_{\tau=-\lambda\xi^2}}\les&
\lambda^{-1/2}(\|u\|^2_{X^{s,b}_{\tau=-\lambda\xi^2}}+\|u'\|^2_{X^{s,b}_{\tau=-\lambda\xi^2}})\|u-u'\|_{X^{s,b}_{\tau=-\lambda\xi^2}}.
\end{align*}
\end{lemma}
\begin{proof}
We just proof the first one when $s=0$ for example. By H\"older
inequality and \eqref{str}, we have
$$
\||u|^2u\|_{X^{0,-a}_{\tau=-\lambda\xi^2}}\leq \||u|^2u\|_{L^2}=
\|u\|^3_{L^6_{t,x}}\leq C
\lambda^{-1/2}\|u\|^3_{X^{0,b}_{\tau=-\lambda\xi^2}}
$$
Thus we finish the proof.
\end{proof}

We also need a bilinear estimates. For $u,v\in \bar{F}^s$ we define the bilinear operator
\begin{align*}
&B(u,v)=\psi({t/4})\int_0^tW(t-\tau)\partial_x\big(\psi^2(\tau)u(\tau)\cdot
v(\tau)\big)d\tau.
\end{align*}
The following lemma is due to the first author \cite{G}, which is
the key to get the global well-posedness for KdV in $H^{-3/4}$.
\begin{lemma}[Proposition 4.2, \cite{G}]\label{nonK}
Assume $-3/4\leq s\leq 0$. Then there exists $C>0$ such that
\begin{eqnarray}\label{eq:bilinearbd}
\norm{B(u,v)}_{\bar{F}^s}\leq
C(\norm{u}_{\bar{F}^s}\norm{v}_{\bar{F}^{-3/4}}+\norm{u}_{\bar{F}^{-3/4}}\norm{v}_{\bar{F}^s})
\end{eqnarray}
hold for any $u,v\in \bar{F}^s$.
\end{lemma}

Now we are ready to prove Theorem \ref{thmlwp}. We consider first \eqref{S-KdV1} under the condition \eqref{indata}.
We may assume $\alpha= \beta= \gamma=1$. From Duhamel's principle, the Cauchy problem \eqref{S-KdV1} is equivalent to the following integral equation system
\begin{eqnarray}\label{inteS-KdV}
\begin{array}{l}
u(t)=U_{\lambda}(t)\phi_1(x)-i\int_{0}^{t}U_\lambda(t-t')[\lambda(uv)(t')+\lambda^{-1}(|u|^2u)(t')]dt'\\
v(t)=V(t)\phi_2(x)+\int_{0}^{t}V(t-t')\Big[-\frac{1}{2}\partial_x(v^2)(t')+\partial_x(|u|^2)(t')\Big]dt'.
\end{array}
\end{eqnarray}
To study the local existence, it suffices to study the following time localized system
\begin{eqnarray}\label{CIe}
\begin{array}{l}
u(t)=\psi_1(t)U_{\lambda}(t)\phi_1(x)-i\psi_T(t)\int_{0}^{t}U_\lambda(t-t')[\lambda(uv)(t')+\lambda^{-1}(|u|^2u)(t')]dt'\\
v(t)=\psi_1(t)V(t)\phi_2(x)-\frac{1}{2}\psi_4(t)\int_{0}^{t}V(t-t')\partial_x(v^2)(t')dt'\\
\quad \quad \quad
+\psi_T(t)\int_{0}^{t}V(t-t')\partial_x(|u|^2)(t')dt'.
\end{array}
\end{eqnarray}
It is easy to see that if $(u,v)$ solve the system \eqref{CIe} for all $t\in \R$, then
$(u,v)$ also solve the system \eqref{inteS-KdV} for $|t|\leq
T$.

We follow the similar argument as the one given in \cite{cl07} to
construct our solution spaces. We consider the following function
space where we seek our solution:
$$
\Sigma_{\theta}:=\Big\{(u,v)\in X^{0,
1/2+\theta}_{\tau=-\lambda\xi^2}\times \bar{F}^{-3/4}_{\tau=\xi^3};
\quad\|u\|_{X^{0, 1/2+\theta}_{\tau=-\lambda\xi^2}}\leq M \text{ and
} \|v\|_{\bar{F}^{-3/4}_{\tau=\xi^3}}\leq \epsilon_0\Big\},
$$
where $0<\theta, \epsilon_0 \ll 1$ and $M>0$, will be chosen later.
Furthermore, $\Sigma_\theta$ is a complete metric space with norm
$$
\|(u,v)\|_{\Sigma_\theta}=\|u\|_{X^{0,
1/2+\theta}_{\tau=-\lambda\xi^2}}+\|v\|_{\bar{F}^{-3/4}_{\tau=\xi^3}}.
$$
For $(u,v)\in \Sigma_\theta$, we define the maps $\Phi=\Phi_1\times\Phi_2$
\begin{align*}
&\Phi_1(u,v)=\psi_1(t)U_{\lambda}(t)\phi_1(x)-i\psi_T(t)\int_{0}^{t}U_\lambda(t-t')[\lambda(uv)(t')+\lambda^{-1}(|u|^2u)(t')]dt'
\end{align*}
and
\begin{align*}
\Phi_2(u,v)=&\psi_1(t)V(t)\phi_2(x)-\frac{1}{2}\psi_4(t)\int_{0}^{t}V(t-t')\partial_x(v^2)(t')dt'\\
&+\psi_T(t)\int_{0}^{t}V(t-t')\partial_x(|u|^2)(t')dt'.
\end{align*}
Then from Lemma \ref{line}-\ref{nonK}, we get
\begin{align*}
\|\Phi_1(u,v)\|_{X_{\tau=-\lambda\xi^2}^{0,1/2+\theta}}\leq&c_0\|\phi_1\|_{L^2}+c_1T^\theta\Big[\lambda\|uv\|_{X_{\tau=-\lambda\xi^2}^{0,-1/2+2\theta}}
 +\lambda^{-1}\||u|^2u\|_{X_{\tau=-\lambda\xi^2}^{0,-1/2+2\theta}}\Big]
\\
\leq&c_0\|\phi_1\|_{L^2}+c_1T^\theta\Big[\lambda\|u\|_{X_{\tau=-\lambda\xi^2}^{0,1/2+\theta}}\|v\|_{\bar{F}_{\tau=\xi^3}^{-3/4}}
 +\lambda^{-3/2}\|u\|^3_{X_{\tau=-\lambda\xi^2}^{0,1/2+\theta}}\Big]
 \\
\leq&c_0\|\phi_1\|_{L^2}+c_1T^\theta\big[\lambda M\epsilon_0
 +\lambda^{-2}M^3\big],
\end{align*}
and
\begin{align*}
\|\Phi_2(u,v)\|_{\bar{F}_{\tau=\xi^3}^{-3/4}}\leq&c_0\|\phi_2\|_{H^{-3/4}}+\frac{1}{2}
     \Big\|\psi(\frac{t}{4})\int_{0}^{t}V(t-t')\partial_x(v^2)(t')dt'\Big\|_{\bar{F}_{\tau=\xi^3}^{-3/4}}\\
&+\Big\|\psi_T(t)\int_{0}^{t}V(t-t')\partial_x(|u|^2)(t')dt'\Big\|_{X_{\tau=\xi^3}^{-3/4,
1/2+\theta}}
\\\leq&c_0\|\phi_2\|_{H^{-3/4}}+c_2\|v\|^2_{\bar{F}_{\tau=\xi^3}^{-3/4}}
 +c_2T^\theta\|\partial_x|u|^2\|_{X_{\tau=\xi^3}^{-3/4, -1/2+2\theta}}
\\
\leq&c_0\|\phi_2\|_{H^{-3/4}}+c_2\|v\|^2_{\bar{F}_{\tau=\xi^3}^{-3/4}}
 +c_2\lambda^{-1/2}T^\theta\|u\|^2_{X_{\tau=-\lambda\xi^2}^{0,1/2+\theta}}
 \\
\leq&c_0\|\phi_2\|_{H^{-3/4}}+c_2\epsilon_0^2
 +c_2\lambda^{-1/2}T^\theta M^2.
\end{align*}
Now taking $M=2c_0 \|\phi_1\|_{L^2}$ and $\epsilon_0\ll 1$ satisfy
$c_2\epsilon_0=1/8$, then we have that
\begin{align*}
\|\Phi_1(u,v)\|_{X_\lambda^{0,1/2+\theta}}\leq&
\frac{M}{2}+c_1T^\theta\big[\lambda M\epsilon_0+\lambda^{-3/2}M^3\big],\\
\|\Phi_2(u,v)\|_{\bar{F}^{-3/4}} \leq&
\frac{\epsilon_0}{2}+c_2\epsilon_0^2
 +c_2\lambda^{-1/2}T^\theta M^2.
\end{align*}
Therefore, we choose $T>0$ such that
\begin{align}\label{exist1}
T^\theta\ll \frac{\lambda^{1/2}}{M^2}\min(\epsilon_0,\lambda),
\end{align}
then $(\Phi_1(u,v),\Phi_2(u,v))\in \Sigma_\theta$. Similarly we have that
\begin{align*}
\|\Phi_1(u,v)-\Phi_1(u',v')\|_{X_{\tau=-\lambda\xi^2}^{0,1/2+\theta}}\leq&
c_3\lambda^{-3/2}T^\theta\big(M+\epsilon_0+M^2\big)\\
&\cdot[\|u-u'\|_{X_{\tau=-\lambda\xi^2}^{0,1/2+\theta}}+\|v-v'\|_{\bar{F}_{\tau=\xi^3}^{-3/4}}],\\
\|\Phi_2(u,v)-\Phi_2(u',v')\|_{\bar{F}_{\tau=\xi^3}^{-3/4}}\leq&
c_3(\lambda^{-1/2}T^\theta M+\epsilon_0)\\
&\cdot[\|u-u'\|_{X_{\tau=-\lambda\xi^2}^{0,1/2+\theta}}+\|v-v'\|_{\bar{F}_{\tau=\xi^3}^{-3/4}}].
\end{align*}
Thus for $T$ in \eqref{exist1} we get
$$
\|(\Phi_1(u,v),\Phi_2(u,v))-(\Phi_1(u',v'),\Phi_2(u',v'))\|_{\Sigma_\theta}\leq
\frac{1}{2}\|(u,v)-(u',v')\|_{\Sigma_\theta}
$$
Therefore the map $\Phi_1\times\Phi_2:\Sigma_\theta\longrightarrow
\Sigma_{\theta}$ is a contraction mapping, and we obtain a unique
fixed point in $\Sigma_\theta$ which solves the equation \eqref{CIe}.

Now we prove Theorem \ref{thmlwp}. Using the scaling \eqref{scaling} by taking
$$\lambda=\Big(\frac{\epsilon_0}{2(\norm{v_0}_{H^{-3/4}}+1)}\Big)^{4/3},$$
then we see
$$M=2c_0 \|\phi_1\|_{L^2}=2c_0 \lambda^{3/2}\|u_0\|_{L^2}, \quad \norm{\phi_2}\leq \epsilon_0.$$
Thus we get local existence of the solution to \eqref{S-KdV1} on $[0,T_1]$ for
\[T_1\sim \brk{\frac{\lambda^{3/2}}{M^2}}^{1/\theta}=\brk{\frac{1}{\lambda^{3/2}\|u_0\|_{L^2}^2}}^{1/\theta}.\]
Therefore, for the original system \eqref{S-KdV}, we get the local existence of the solution on $[0,T]$ for
\[T=\lambda^3T_1\sim T(\norm{u_0}_{L^2},\norm{v_0}_{H^{-3/4}}).\]
Therefore, we obtain existence. By standard arguments, we can prove
the rest parts of Theorem \ref{thmlwp}.

\end{document}